# TRANSLATED POISSON APPROXIMATION USING EXCHANGEABLE PAIR COUPLINGS

By Adrian Röllin[1]

*University of Zurich*

It is shown that the method of exchangeable pairs introduced by Stein [*Approximate Computation of Expectations* (1986) IMS, Hayward, CA] for normal approximation can effectively be used for translated Poisson approximation. Introducing an additional smoothness condition, one can obtain approximation results in total variation and also in a local limit metric. The result is applied, in particular, to the anti-voter model on finite graphs as analyzed by Rinott and Rotar [*Ann. Appl. Probab.* **7** (1997) 1080–1105], obtaining the same rate of convergence, but now for a stronger metric.

**1. Introduction.** Let $W$ be a random variable with

$$(1.1) \qquad \mathbb{E}W = \mu \quad \text{and} \quad \operatorname{Var}W = \sigma^2 < \infty.$$

Stein [17] introduced a method (which is commonly called the *exchangeable pairs* approach) to approximate $W_c := (W - \mu)/\sigma$ by the standard normal distribution; Rinott and Rotar [14] then generalized the result and successfully applied it to weighted $U$-statistics and the antivoter model. Their results imply convergence to the standard normal distribution in the Kolmogorov and even in some stronger metrics; however, in this context, they do not provide approximations in the total variation metric or prove local limit-like results.

We will consider such results in this paper in the special case, in which $W$ is integer valued, the most common situation being the one where $W$ is a sum of random indicators. As the total variation distance between $W$ and the normal distribution will always be 1, we will instead use a translated Poisson distribution as approximation, having the same support as $W$ and

---

Received August 2006; revised March 2007.

[1]Supported in part by Schweizerischer Nationalfondsprojekt 20-107935/1.

matching the first two moments of $W$ as well as possible. If not otherwise stated, we will assume throughout that $\sigma^2 \to \infty$. This actually implies that the approximating probability distributions are not converging to a limiting distribution, but the accuracy of our approximations nonetheless increases as $\sigma^2$ becomes large. The total variation metric is invariant under scaling, so that working with $W_c$ would bring no benefit. Besides the total variation metric, we will also consider a metric from which local limit approximations can be obtained.

We note that in the framework of Stein's method there are other approaches to replacing the normal distribution by discrete analogues. In [12] a distribution with support on $\mathbb{Z}$ is constructed, with the advantage of having no truncation and rounding effect but at the cost of a somewhat more complicated Stein operator. There, approximation theorems are provided using the so-called *zero biasing* approach introduced in [9]. In [16] a translated binomial distribution is used, and in [4] a special translated signed compound Poisson distribution, both in the context of the so-called *local approach.* In the discrete setting, exchangeable pairs have also been successfully used in [8] for Poisson approximation in total variation.

The rest of the paper is organized as follows. In the next section we recall the setup of the exchangeable pairs approach in the context of normal approximation. We then introduce a simple smoothing condition under which it is possible to obtain the stronger total variation bounds for translated Poisson approximation. In Section 3 we state and prove the main approximation theorem, Theorem 3.1, which is the discrete equivalent to Theorem 1.2 of [14]. We also prove a second general result, Theorem 3.11, from which, under an additional assumption, more accurate rates can be obtained for local limit results. In Section 4 some applications are given, among others to the anti-voter model.

## 2. Exchangeable pairs for normal approximation and a smoothness condition.
We call a pair of random variables $(W, W')$ exchangeable if $\mathscr{L}(W, W') = \mathscr{L}(W', W)$. As in [17] and [14], assume now that there is a positive number $\lambda < 1$ and a random variable $R$ such that

$$(2.1) \qquad \mathbb{E}^W(W' - \mu) = (1 - \lambda)(W - \mu) + R,$$

holds, where $\mathbb{E}^W$ denotes the conditional expectation with respect to $W$. Of course, one can always find $R$ to satisfy (2.1), so $R$ must be thought of as being small for the approximation to be successful. Note that (2.1) implies $\mathbb{E}R = 0$.

If the pair $(W, W')$ can be chosen such that condition (2.1) holds and $\mathbb{E}^W(W' - W)^2$ does not fluctuate too much, convergence of $W_c$ to the standard normal distribution will follow in the Kolmogorov metric. As the behavior of the difference $W' - W$ is mainly responsible for the quality of the approximation, it is an obvious starting point to introduce a smoothness condition, to make sure that the local perturbations of $W$ are not too strong.

Rinott and Rotar [14] propose to choose $W$ and $W'$ as two successive steps of a reversible Markov chain with stationary distribution $\mathscr{L}(W)$. Then, condition (2.1) states that a particle on $\mathbb{Z}$ obeying the transition rules of such a Markov chain is forced to have an (almost) linear drift to the center. Now $\mathbb{E}^{W=k}(W'-W)^2$ is the average of the squared jump size of the Markov chain if the particle is in $k$, so that, for a good normal approximation, the average jump size of the particle must not fluctuate too much with varying $k$. It is clear that, under these conditions, the particle may still behave irregularly on a local scale, for instance, the particle could still make only jumps of size two and thus stay on the odd or even integers, such that an approximation with a distribution on $\mathbb{Z}$ with span 1 will not be successful in total variation.

Thus, in addition to (2.1), we assume further that

$$(2.2) \qquad\qquad W' - W \in \{-1, 0, +1\},$$

and we will see that this seems to be an appropriate condition. Note that under condition (2.2) the corresponding Markov chain does not need to be reversible for $(W, W')$ to be a exchangeable pair; see Lemma 1.1 of [14].

Condition (2.2) is in sharp contrast to other approaches using Stein's method for the translated Poisson distribution such as [7, 15] or [3], where an embedded sum of independent random variables within $W$ is used for an explicit smoothing argument; in contrast, the smoothing effect of (2.2) will enter only implicitly into the proof of the main result. As we are restricted to the integers, we cannot arbitrarily shift a Poisson distribution with a given variance to fit the mean, so some care is needed here. We say that an integer valued random variable $Y$ has a *translated Poisson distribution* with parameters $\mu$ and $\sigma^2$ and write

$$\mathscr{L}(Y) = \mathrm{TP}(\mu, \sigma^2)$$

if $\mathscr{L}(Y - \mu + \sigma^2 + \gamma) = \mathrm{Po}(\sigma^2 + \gamma)$, where $\gamma = \langle \mu - \sigma^2 \rangle$ and $\langle x \rangle = x - \lfloor x \rfloor$ denotes the fractional part of $x$; in particular, $\mathrm{TP}(\sigma^2, \sigma^2) = \mathrm{Po}(\sigma^2)$. So, approximating $W$ with $\mathrm{TP}(\mu, \sigma^2)$, we can fit the mean exactly, but note that for the variance we have $\sigma^2 \leq \mathrm{Var}\, Y = \sigma^2 + \gamma \leq \sigma^2 + 1$. This will, however, cause no further problems as the order of error of this mismatch is $O(\sigma^{-2})$; see also Remark 3.5 below.

Throughout the paper, we shall be concerned with two metrics for probability distributions, the total variation metric $d_{\mathrm{TV}}$ and the local limit metric $d_{\mathrm{loc}}$, where, for two probability distributions $P$ and $Q$ given by the point probabilities $\{p_k, k \in \mathbb{Z}\}$ and $\{q_k, k \in \mathbb{Z}\}$ respectively,

$$d_{\mathrm{TV}}(P, Q) := \sup_{A \subset \mathbb{Z}} |P(A) - Q(A)| = \tfrac{1}{2} \sum_{k \in \mathbb{Z}} |p_k - q_k|,$$

$$d_{\mathrm{loc}}(P, Q) := \sup_{k \in \mathbb{Z}} |p_k - q_k|.$$

## 3. Main results.

THEOREM 3.1. *Assume that $(W, W')$ is an exchangeable pair with values on the integers and which satisfies (1.1), (2.1) and (2.2). Then, with $S = S(W) = \mathbb{P}[W' = W + 1|W]$ and $q_{\max} = \max_{k \in \mathbb{Z}} \mathbb{P}[W = k]$,*

$$(3.1) \qquad d_{\mathrm{TV}}(\mathscr{L}(W), \mathrm{TP}(\mu, \sigma^2)) \leq \frac{\sqrt{\mathrm{Var}\, S}}{\lambda \sigma^2} + \frac{2\sqrt{\mathrm{Var}\, R}}{\lambda \sigma} + \frac{2}{\sigma^2},$$

$$(3.2) \qquad d_{\mathrm{loc}}(\mathscr{L}(W), \mathrm{TP}(\mu, \sigma^2)) \leq \frac{2\sqrt{q_{\max} \mathrm{Var}\, S}}{\lambda \sigma^2} + \frac{2 q_{\max} \sqrt{\mathrm{Var}\, R}}{\lambda \sigma}$$
$$+ \frac{\sqrt{\mathrm{Var}\, R}}{\lambda \sigma^2} + \frac{2}{\sigma^2}.$$

REMARK 3.2. *In some of the applications, instead of $S(W) = \mathbb{P}[W' = W + 1|W]$, we will estimate the variance of a more general random variable $S^* = S^*(X) := \mathbb{P}[W' = W + 1|X]$ for some random variable $X$ such that the corresponding $\sigma$-algebras satisfy $\sigma(W) \subset \sigma(X)$ and then use the basic fact that $\mathrm{Var}\, S \leq \mathrm{Var}\, S^*$.*

EXAMPLE 3.3. *To illustrate the above theorem, we apply it to the Poisson-binomial distribution. To this end let $J = (J_1, \ldots, J_n)$ be a sequence of independent indicators with $\mathbb{E} J_i = p_i$ and $W = \sum_{i=1}^n J_i$, thus, $\mu = \sum_{i=1}^n p_i$ and $\sigma^2 = \sum_{i=1}^n p_i(1 - p_i)$. We use the standard construction of [17] to obtain an exchangeable pair. Let $J_1^*, \ldots, J_n^*$ be independent copies of the $J_i$ and let $K$ be uniformly distributed over $\{1, \ldots, n\}$. Then, with $W' = W - J_K + J_K^*$, it is easy to check that $(W', W)$ is an exchangeable pair, satisfying (2.1) with $R \equiv 0$ and $\lambda = 1/n$ and, clearly, (2.2) is also satisfied. So,*

$$(3.3) \qquad \begin{aligned} S^*(J) &:= \mathbb{E}^J I[W' - W = 1] \\ &= \frac{1}{n} \sum_{i=1}^n \mathbb{E}^J I[J_i = 0, J_i^* = 1] \\ &= \frac{1}{n} \sum_{i=1}^n (1 - J_i) \mathbb{E}^J J_i^* \\ &= \frac{1}{n} \sum_{i=1}^n (1 - J_i) p_i. \end{aligned}$$

*Thus, $\mathrm{Var}\, S^* = n^{-2} \sum_{i=1}^n p_i^3 (1 - p_i)$, and, by Remark 3.2, (3.1) yields*

$$(3.4) \qquad d_{\mathrm{TV}}(\mathscr{L}(W), \mathrm{TP}(\mu, \sigma^2)) \leq \frac{2 + \sqrt{\sum p_i^3 (1 - p_i)}}{\sum p_i (1 - p_i)}.$$

Assume now that the $p_i$ are bounded away from 0 and 1, so that $\sigma^2 \asymp n$ as $n \to \infty$. Then (3.4) is of the correct order $O(n^{-1/2})$. This also implies that $q_{max} = O(n^{-1/2})$ (see Corollary 3.9 below) so that (3.2) yields $d_{loc}(\mathscr{L}(W), \text{TP}(\mu, \sigma^2)) = O(n^{-3/2})$, which in contrast is not optimal. We will improve this bound using Theorem 3.11 below.

In the above case of the Poisson-binomial distribution, Corollary 2.1 of [7] seems to be better in constant than (3.1). For instance, for the binomial distribution, we have

$$d_{\text{TV}}(\text{Bi}(n, p), \text{TP}(\mu, \sigma^2)) \leq C\sqrt{\frac{p}{n(1-p)}} + \frac{2}{np(1-p)},$$

where [7] obtain $C = 0.93$ and (3.1) yields $C = 1$.

REMARK 3.4. Theorem 3.1 is a direct analogue of Theorem 1.2 of [14]. However, the first term in (3.1) is slightly different in quality from Theorem 1.2 of [14], as can be seen by comparing the result of their Theorem 1.3 for the anti-voter model with estimate (4.4) below.

REMARK 3.5. The additional $2/\sigma^2$ in (3.1) and (3.2) appears because the Poisson distribution cannot take negative values, and because the translation must be integer valued. Depending on the problem at hand, this error term can be further reduced or even be omitted by replacing estimates (3.12) and (3.17) in the proof below. For example, to obtain the best possible total variation estimates from (3.1) in the Poisson-binomial case, recall from Section 2 that $\gamma = \langle \mu - \sigma^2 \rangle = \langle \sum p_i^2 \rangle$ and $s = \lfloor \mu - \sigma^2 \rfloor = \mu - \sigma^2 - \gamma = \sum p_i^2 - \langle \sum p_i^2 \rangle$. From (2.8) of [7] we obtain for (3.12)

$$\mathbb{P}[W < s] \begin{cases} \leq e^{-\sigma^2/4}, & \text{if } s > 0, \\ = 0, & \text{if } s = 0. \end{cases}$$

For the last term in (3.16), we have

$$|\mathbb{E}\gamma\Delta\tilde{g}(W)| \leq \|\Delta\tilde{g}\|\langle\sum p_i^2\rangle.$$

Using the first inequality of the estimate of $\|\Delta\tilde{g}\|$ in (3.13) and applying this in the above estimate and also in (3.23), we obtain

$$
\begin{aligned}
(3.5) \quad & d_{\text{TV}}(\mathscr{L}(W), \text{TP}(\mu, \sigma^2)) \\
& \leq \frac{1 - e^{-\sigma^2 - \langle\sum p_i^2\rangle}}{\sigma^2 + \langle\sum p_i^2\rangle}\left(\sqrt{\sum p_i^3(1-p_i)} + \langle\sum p_i^2\rangle\right) \\
& \quad + I\left[\sum p_i^2 \geq 1\right]e^{-\sigma^2/4}.
\end{aligned}
$$

This estimate now covers also the regime of Poisson approximation. However, (3.5) is larger in constant than previous results and one would have to go back to the proof of the theorem to reproduce the bounds of [1] and [8]; see Remark 3.8.

REMARK 3.6. As becomes clear from equation (3.22) in the proof of Theorem 3.1, there is a close connection between the random variable $S = S(W)$ and the so-called $w$-functions as examined, for example, in [5] and [6] for the normal and the Poisson distributions. In the case of the standard normal distribution, their problem is as follows: for a given random variable $X$ with $\mathbb{E}X = 0$ and $\operatorname{Var}X = 1$, find a function $w : \mathbb{R} \to \mathbb{R}$ such that

$$(3.6) \qquad \mathbb{E}\{Xf(X)\} = \mathbb{E}\{w(X)f'(X)\}$$

holds for a large set of functions $f$. For the translated Poisson distribution, the corresponding equation is

$$(3.7) \qquad \mathbb{E}\{(W - \mu)f(W)\} = \mathbb{E}\{w(W)\Delta f(W)\},$$

and it is indeed satisfied for any $W$ as in Theorem 3.1 if $R \equiv 0$ and if we choose $w(W) = S(W)/\lambda$. Unfortunately, it is often difficult to give an explicit expression for $S$ as a function of $W$. However, if we allow $w(W)$ in (3.7) to be replaced by a more general random variable, we see from (3.22) that we can use the random variable $S^*(X)/\lambda$ from Remark 3.2 instead. For instance, for the anti-voter model as discussed in the next section, $S^*(X)$ has the nice and explicit representation (4.7).

Instead of (3.6), one can also formulate the problem of finding a random variable $X^z$ such that

$$(3.8) \qquad \mathbb{E}\{Xf(X)\} = \mathbb{E}f'(X^z),$$

which leads to the so-called *zero biasing* approach. There is a close connection between this and the exchangeable pairs approach; see [10] and references therein, and for more general versions of (3.8), see [11].

Before proving Theorem 3.1, we give a short introduction into Stein's method for translated Poisson approximation. The starting point is the Stein–Chen method for the Poisson distribution as presented in detail by Barbour, Holst and Janson [2].

Let $W$ satisfy (1.1) and let $s = \lfloor \mu - \sigma^2 \rfloor$ and $\gamma = \langle \mu - \sigma^2 \rangle$, where $\langle x \rangle = x - \lfloor x \rfloor$ denotes the fractional part of $x$. Note that, if $Y \sim \mathrm{TP}(\mu, \sigma^2)$, then $Y - s \sim \mathrm{Po}(\sigma^2 + \gamma)$. Let $\mathcal{A}g(j) = (\sigma^2 + \gamma)g(j+1) - jg(j)$ be the usual Stein operator for the Poisson distribution with mean $\sigma^2 + \gamma$, and for $A \subset \mathbb{Z}_+ := \{0, 1, 2, \ldots\}$, let $g_A : \mathbb{Z} \to \mathbb{R}$ be the solution of the following:

$$(3.9) \qquad \text{(i)} \quad g(j) = 0 \qquad \text{for all } j \leq 0,$$

$$(3.10) \qquad \text{(ii)} \quad \mathcal{A}g(j) = I[j \in A] - \mathrm{Po}(\sigma^2 + \gamma)\{A\} \qquad \text{for all } j \geq 0.$$

We can thus bound the total variation distance as

(3.11)
$$
\begin{aligned}
d_{\mathrm{TV}}(&\mathscr{L}(W), \mathrm{TP}(\mu, \sigma^2)) \\
&= d_{\mathrm{TV}}(\mathscr{L}(W-s), \mathrm{Po}(\sigma^2+\gamma)) \\
&= \sup_{B \subset \mathbb{Z}} |\mathbb{E} I[W-s \in B] - \mathrm{Po}(\sigma^2+\gamma)\{B\}| \\
&\leq \sup_{A \subset \mathbb{Z}_+} |\mathbb{E} \mathcal{A} g_A(W-s)| + \mathbb{P}[W-s < 0].
\end{aligned}
$$

The last term in (3.11) can be bounded using Chebyshev's inequality as

(3.12)
$$
\begin{aligned}
\mathbb{P}[W-s < 0] &= \mathbb{P}[W-\mu < -(\sigma^2+\gamma)] \\
&\leq \mathbb{P}[|W-\mu| > \sigma^2+\gamma] \leq \frac{1}{\sigma^2}.
\end{aligned}
$$

From [2], Lemma 1.1.1, we have the well-known bounds on the supremum norm of $g_A$,

(3.13)
$$
\|g_A\| \leq (\sigma^2+\gamma)^{-1/2} \leq \sigma^{-1},
$$
$$
\|\Delta g_A\| \leq \frac{1-e^{-\sigma^2-\gamma}}{\sigma^2+\gamma} \leq \sigma^{-2},
$$

where $\Delta g_A(j) := g_A(j+1) - g_A(j)$. If $A = \{k\}$ for some $k \in \mathbb{Z}$, we even have

(3.14)
$$
\|g_{\{k\}}\| \leq \sigma^{-2}.
$$

For the proof of the results in the $d_{\mathrm{loc}}$ metric, we will also need the following nonstandard but simple result.

LEMMA 3.7. *Let $g_i$ be the solution of* (3.9)–(3.10) *for $A = \{i\}$. Then*

(3.15)
$$
\sum_k |\Delta g_i(k)| \leq 2\sigma^{-2}, \qquad \sum_k (\Delta g_i(k))^2 \leq 4\sigma^{-4}.
$$

PROOF. Recall from [2], proof of Lemma 1.1.1, that $g_i(k)$ is negative and decreasing in $0 \leq k \leq i$ and positive and decreasing in $k > i$ with the only positive jump in $i$ satisfying

$$
|\Delta g_i(i)| \leq (\sigma^2+\gamma)^{-1} \leq \sigma^{-2}.
$$

From this, it is easy to see that the first bound of (3.15) holds and the second bound is then immediate. □

With $\tilde{g}_A(j) := g_A(j-s)$, we can rewrite the Stein operator $\mathcal{A}$ as

(3.16)
$$
\begin{aligned}
\mathcal{A} g_A(W-s) &= (\sigma^2+\gamma) g_A(W-s+1) - (W-s) g_A(W-s) \\
&= \sigma^2 \Delta \tilde{g}_A(W) - (W-\mu) \tilde{g}_A(W) + \gamma \Delta \tilde{g}_A(W).
\end{aligned}
$$

The bounds on $\tilde{g}_A$ are of course the same as on $g_A$ in (3.13)–(3.15). Thus, the expectation of the last term in (3.16) is easily bounded by

$$(3.17) \qquad |\mathbb{E}\{\gamma\Delta\tilde{g}_A(W)\}| \leq \gamma\sigma^{-2} \leq \sigma^{-2}.$$

Inserting (3.16) into (3.11) and invoking the bounds (3.12) and (3.17), we obtain

$$
\begin{aligned}
(3.18) \qquad & d_{\mathrm{TV}}(\mathscr{L}(W), \mathrm{TP}(\mu, \sigma^2)) \\
& \leq \sup_{A \subset \mathbb{Z}_+} |\mathbb{E}\{\sigma^2\Delta\tilde{g}_A(W) - (W-\mu)\tilde{g}_A(W)\}| + 2\sigma^{-2};
\end{aligned}
$$

the same estimate holds for $d_{\mathrm{loc}}$ but with the supremum taken only over the sets $A = \{i\}$ for $i \in \mathbb{Z}_+$.

PROOF OF THEOREM 3.1. We only have to bound the supremum in (3.18). In [17] it was shown that, if $F$ satisfies $F(w, w') = -F(w', w)$ for all $w$ and $w'$, exchangeability implies $\mathbb{E}F(W, W') = 0$. Define the random variable $D := W' - W$ and the function $F(w, w') := (w' - w)(g(w') + g(w))$ for $g \equiv \tilde{g}_A$ and note that, from (2.1), $\mathbb{E}^W D = -\lambda(W-\mu) + R$. This yields

$$
\begin{aligned}
(3.19) \qquad 0 &= \mathbb{E}F(W, W') = \mathbb{E}\{D(2g(W) + g(W') - g(W))\} \\
&= -2\lambda\mathbb{E}\{(W-\mu)g(W)\} + 2\mathbb{E}\{Rg(W)\} + \mathbb{E}\{D(g(W') - g(W))\}.
\end{aligned}
$$

Note now that, for $D_i := I[D = i]$, $i \in \{-1, +1\}$, we can write

$$D(g(W') - g(W)) = D_{+1}\Delta g(W) + D_{-1}\Delta g(W-1),$$

and further, using exchangeability,

$$
\begin{aligned}
(3.20) \qquad \mathbb{E}\{D_{-1}\Delta g(W-1)\} &= \mathbb{E}\{I[W' - W = -1]\Delta g(W-1)\} \\
&= \mathbb{E}\{I[W - W' = 1]\Delta g(W')\} \\
&= \mathbb{E}\{D_{+1}\Delta g(W)\},
\end{aligned}
$$

thus,

$$(3.21) \qquad \mathbb{E}\{D(g(W') - g(W))\} = 2\mathbb{E}\{D_{+1}\Delta g(W)\}.$$

Together with (3.19) this yields

$$(3.22) \qquad \mathbb{E}\{(W-\mu)g(W)\} = \frac{\mathbb{E}\{D_{+1}\Delta g(W)\}}{\lambda} + \frac{\mathbb{E}\{Rg(W)\}}{\lambda}.$$

Note now that, by exchangeability, $\mathbb{E}D_{+1} = \mathbb{E}D_{-1}$ and, hence, that

$$
\begin{aligned}
(3.23) \qquad \mathbb{E}D_{+1} &= \tfrac{1}{2}\mathbb{E}(W' - W)^2 \\
&= \tfrac{1}{2}[\mathbb{E}(W' - \mu)^2 - 2\mathbb{E}\{(W' - \mu)(W - \mu)\} + b\mathbb{E}(W - \mu)^2] \\
&= \lambda\sigma^2 + \mathbb{E}\{(W-\mu)R\} =: \lambda\sigma^2 + a,
\end{aligned}
$$

from (2.1); then use (3.22) to express the expectation in (3.18) as

$$\mathbb{E}\{(W - \mu)g(W) - \sigma^2 \Delta g(W)\}$$

$$= \mathbb{E}\{(W - \mu)g(W) - (\sigma^2 + \lambda^{-1}a)\Delta g(W)\} + \lambda^{-1}a\mathbb{E}\Delta g(W)$$

$$= \mathbb{E}\{(D_{+1}\lambda^{-1} - \sigma^2 - \lambda^{-1}a)\Delta g(W)\} + \lambda^{-1}\mathbb{E}\{Rg(W)\} + \lambda^{-1}a\mathbb{E}\Delta g(W)$$

$$=: B_1 + B_2 + B_3.$$

Now, recall that $S = \mathbb{E}^W D_{+1}$, and thus, with the estimates

$$(3.24) \qquad |B_1| \leq \|\Delta g\|\lambda^{-1}\mathbb{E}|S - \mathbb{E}S| \leq \|\Delta g\|\lambda^{-1}\sqrt{\operatorname{Var} S},$$

$$(3.25) \qquad |B_2| \leq \|g\|\lambda^{-1}\mathbb{E}|R| \leq \|g\|\lambda^{-1}\sqrt{\operatorname{Var} R},$$

$$(3.26) \qquad |B_3| \leq \|\Delta g\|\lambda^{-1}\mathbb{E}|(W - \mu)R| \leq \|\Delta g\|\lambda^{-1}\sigma\sqrt{\operatorname{Var} R},$$

and the bounds (3.13), (3.1) follows.

To prove (3.2), we also use (3.18), but now we take the supremum only over all subsets $A = \{i\}$ for $i \in \mathbb{Z}$. Writing $g \equiv \tilde{g}_{\{i\}}$ and following the proof as for $d_{\mathrm{TV}}$ above, the bound on (3.25) remains and recalling (3.14), the third term in (3.2) follows. We thus need only refine the bounds on $B_1$ and $B_3$. Note that by the Cauchy–Schwarz inequality

$$|B_1| \leq \lambda^{-1}\sqrt{\operatorname{Var} S}\sqrt{\mathbb{E}(\Delta g(W))^2}.$$

Using Lemma 3.7, the latter expectation can be bounded by

$$(3.27) \qquad \begin{aligned} \mathbb{E}(\Delta g(W))^2 &= \sum_k (\Delta g(k))^2 \mathbb{P}[W = k] \\ &\leq q_{\max} \sum_k (\Delta g(k))^2 \leq 4\sigma^{-4}q_{\max}, \end{aligned}$$

which implies the first term in (3.2). Using a similar argument on $B_3$, we obtain

$$|B_3| \leq \lambda^{-1}\sigma\sqrt{\operatorname{Var} R}\, q_{\max}\sum_k |\Delta g(k)|,$$

which, together with Lemma 3.7, yields the second term in (3.2).  $\square$

REMARK 3.8.   It is interesting to compare our approach to the one used in [8], who also use exchangeable pairs but for Poisson approximation. As we have $\mathrm{TP}(\mu, \mu) = \mathrm{Po}(\mu)$, it should be expected that we can reproduce their results. This is indeed the case.

Now, assume our conditions (2.1) and (2.2) and assume that we are in the regime of Poisson approximation, that is, $\sigma^2 \approx \mu$. We also assume for the sake of simplicity that $R \equiv 0$. Taking $\mathrm{TP}(\mu, \mu)$ as the approximating distribution, it is easy to see that the Stein operator (3.16) reduces to the classical Stein operator

$$\mathcal{A}g(w) := \mu g(w+1) - wg(w)$$

for $Po(\mu)$ also used in [8]. Using the anti-symmetric function from the proof of Theorem 3.1, we have

$$
\begin{aligned}
(3.28) \qquad 0 &= \mathbb{E}\{(D_{+1} - D_{-1})(g(W') + g(W))\} \\
&= \mathbb{E}\{D_{+1}g(W+1) - D_{-1}g(W)\},
\end{aligned}
$$

where for the second equality we exploited the same argument as in (3.20) and the fact that $W' = W + 1$ if $D_{+1} = 1$. Note that (3.28) is the same equality as in [8], obtained through a different anti-symmetric function. Multiplying (3.28) by an arbitrary constant $c$, we obtain the bound

$$
\begin{aligned}
d_{\mathrm{TV}}(\mathscr{L}(W), Po(\mu)) = d_{\mathrm{TV}}(\mathscr{L}(W), \mathrm{TP}(\mu,\mu)) &= \sup_g |\mathbb{E}\mathcal{A}g(W)| \\
(3.29) \qquad &\leq \sup_g |\mathbb{E}\{(\mu - c\mathbb{E}^W D_{+1})g(W+1) \\
&\qquad\qquad + (W - c\mathbb{E}^W D_{-1})g(W)\}|,
\end{aligned}
$$

where the supremum ranges over the same functions $g$ as in (3.11). Note that, from (3.23), we have $\mathbb{E}^W D_{+1} = \lambda\sigma^2 \approx \lambda\mu$ and this, in conjunction with (2.1), also implies $\mathbb{E}^W D_{-1} \approx \lambda W$, so that, with the choice $c = \lambda^{-1}$, (3.29) is expected to be reasonably small in the regime of Poisson approximation.

Instead of (2.1) and (2.2), in [8] it is only assumed that $(W, W')$ is exchangeable. With this assumption, they prove the same bound (3.29), where again $c$ can be chosen arbitrarily. It is noteworthy that, although differences of $|W' - W|$ larger than 1 are allowed in their approach, again only jumps of size 1 appear in (3.29); this is a consequence of exchangeability.

So, we are able to reproduce the estimates of [8] under our assumptions, by taking $\mathrm{TP}(\mu,\mu)$ instead of $\mathrm{TP}(\mu,\sigma^2)$ as the approximation. However, we have the extra flexibility of being able to match mean and variance separately, so that our approach also works when $\sigma^2$ is not near $\mu$; for instance, if $\sum p_i^2$ is not small in the Poisson-binomial case. In contrast, in [8] they do not assume (2.1), and allow for differences larger than 1.

Using (3.1) with the following corollary, one easily obtains a bound for $q_{\max}$.

COROLLARY 3.9. *For any $\mathbb{Z}$-valued random variable $W$,*

$$q_{\max} \leq d_{\mathrm{TV}}(\mathscr{L}(W), \mathrm{TP}(\mu,\sigma^2)) + \frac{1}{2.3\sigma}.$$

Proof. Just apply Proposition A.2.7 of [2]. □

Remark 3.10. Estimate (3.2) in combination with Corollary 3.9 is enough to obtain a local limit theorem in the applications of the next section. Although it can be easily calculated in many circumstances, the example of the Poisson-binomial distribution shows that the bound on $d_{loc}$ need not be optimal; estimate (3.2) is of order $O(n^{-3/4})$ in the special case of the binomial distribution, in contrast to the true order $O(n^{-1})$. Under additional assumptions on $S$, however, the bound (3.2) can be used to derive the better $d_{loc}$-bound, given in the following theorem. This bound is used in the examples of the Sections 4.1 and 4.2 to obtain the correct order $O(n^{-1})$ of approximation.

Theorem 3.11. Assume the conditions of Theorem 3.1; assume, in addition, that $S$, as a function of $W$, can be extended on $\mathbb{R}$ such that it is Lipschitz continuous. Then,

(3.30)
$$
\begin{aligned}
&d_{loc}(\mathscr{L}(W), \mathrm{TP}(\mu, \sigma^2)) \\
&\quad \leq \frac{2L_S(\sigma^{-3}\mathbb{E}|W - \mu|^3 \vee (d\sigma^{3/2} + 1))}{\lambda\sigma^2} + \frac{2L_S q_{max}}{\lambda\sigma} \\
&\quad + \frac{2q_{max}\sqrt{\mathrm{Var}\,R}}{\lambda\sigma} + \frac{\sqrt{\mathrm{Var}\,R}}{\lambda\sigma^2} + \frac{2}{\sigma^2},
\end{aligned}
$$

where $d$ is the $d_{loc}$-bound (3.2) and $L_S$ is the Lipschitz constant of $S$.

To obtain useful bounds from the above theorem, it is essential that one has a good bound on $L_S$. In the Sections 4.1 and 4.2 and in the special case of the anti-voter model on the complete graph (Example 4.7), this is easily done, because there we know $S$ explicitly. Recall also Example 3.3 for the binomial distribution, that is, $p_i = p$ for some fixed $p$. Then (3.3) yields $S^*(J) = \lambda p(n - W) = S(W)$. Clearly, $L_S = \lambda p$, so that from (3.30) we obtain the correct order $O(n^{-1})$ for the $d_{loc}$-metric. For the general Poisson-binomial and anti-voter models from Section 4.3, however, we only know a more general function $S^*(J)$ with $S(W) = \mathbb{E}^W S^*(J)$ (see Remark 3.2), and it is unclear how to obtain useful bounds on $L_S$ in these cases.

To prove Theorem 3.11, we need the following lemma.

Lemma 3.12. For any $\mu$ and $\sigma^2$, the bound

$$
\mathrm{TP}(\mu, \sigma^2)\{k\}|k - \mu| \leq 1
$$

holds for all $k \in \mathbb{Z}$.

PROOF. Recall from (3.10) that, if $Z \sim \mathrm{TP}(\mu, \sigma^2)$,

$$(3.31) \qquad \mathbb{E}\{(Z - \mu)g(Z) - (\sigma^2 + \gamma)\Delta g(Z)\} = 0$$

for any $g$ for which the expectations exist. With $\pi_k = \mathrm{TP}(\mu, \sigma^2)\{k\}$ and putting $g(\cdot) = I[\cdot = k]$ we obtain from (3.31) the bound

$$
\begin{aligned}
\pi_k |k - \mu| &\leq (\sigma^2 + \gamma)|\pi_{k-1} - \pi_k| \\
&\leq (\sigma^2 + \gamma)d_{\mathrm{loc}}(\mathrm{TP}(\mu + 1, \sigma^2), \mathrm{TP}(\mu, \sigma^2)) \\
&= (\sigma^2 + \gamma)d_{\mathrm{loc}}(\mathscr{L}(Y + 1), \mathscr{L}(Y)),
\end{aligned}
$$

where $Y \sim \mathrm{Po}(\sigma^2 + \gamma)$. The later $d_{\mathrm{loc}}$-distance can easily be bounded using Stein's method for the Poisson distribution, that is, (3.10) in connection with the bound (3.14), which yields $d_{\mathrm{loc}}(\mathscr{L}(Y + 1), \mathscr{L}(Y)) \leq (\sigma^2 + \gamma)^{-1}$. □

PROOF OF THEOREM 3.11. Follow the proof of Theorem 3.1 for the $d_{\mathrm{loc}}$ metric up to the bounds on the $B_i$. The bounds on $|B_2|$ and $|B_3|$ remain. Recalling that $S$ is a function defined on all $\mathbb{R}$, write now $B_1$ as

$$
\begin{aligned}
B_1 &= \lambda^{-1}\mathbb{E}\{(S(W) - \mathbb{E}S(W))\Delta g(W)\} \\
&= \lambda^{-1}\mathbb{E}\{(S(W) - S(\mu))\Delta g(W)\} + \lambda^{-1}\mathbb{E}\{(S(\mu) - S(W))\}\mathbb{E}\Delta g(W) \\
&=: B_{1,1} + B_{1,2}.
\end{aligned}
$$

Exploiting Lipschitz continuity of $S$ and recalling (3.15), we obtain with $q_k = \mathbb{P}[W = k]$

$$|B_{1,2}| \leq \lambda^{-1}\sigma L_S \sum_k q_k |\Delta g(k)| \leq \frac{2L_S q_{\max}}{\lambda \sigma},$$

which is the second term in (3.30). For $B_{1,1}$, we have

$$(3.32) \qquad
\begin{aligned}
|B_{1,1}| &\leq \lambda^{-1} \sum_k q_k |S(k) - S(\mu)||\Delta g(k)| \\
&\leq \lambda^{-1}L_S \sum_k q_k |k - \mu||\Delta g(k)|.
\end{aligned}
$$

We now bound $q_k |k - \mu|$. Assume first that $|k - \mu| > \sigma^{3/2}$; then, by Chebyshev's inequality,

$$q_k \leq \mathbb{P}[W \geq k] \leq \frac{\mathbb{E}|W - \mu|^3}{|k - \mu|^3}\mathbb{P}[|W - \mu| \geq |k - \mu|]$$

and, thus,

$$q_k |k - \mu| \leq \sigma^{-3}\mathbb{E}|W - \mu|^3.$$

On the other hand, if $|k - \mu| \leq \sigma^{3/2}$, observe that

$$q_k \leq d + \mathrm{TP}(\mu, \sigma^2)\{k\}$$

and, hence, using Lemma 3.12,

$$q_k|k - \mu| \leq d\sigma^{3/2} + 1.$$

Thus, (3.32) can be further bounded to

$$|B_{1,1}| \leq \lambda^{-1} L_S(\sigma^{-3}\mathbb{E}|W - \mu|^3 \vee (d\sigma^{3/2} + 1)) \sum_k |\Delta g(k)|$$

and applying again (3.15), this yields the first term in (3.30).  $\square$

The following lemma can be used to estimate the second and third moments of $W$.

LEMMA 3.13.  *Under the assumptions of Theorem* 3.1 *and with* $A = \{w : \mathbb{P}[W = w] > 0\}$ *and* $a := \mathbb{E}\{R(W - \mu)\}$,

$$\lambda^{-1}\left(\inf_{w \in A} S(w) - a\right) \leq \sigma^2 \leq \lambda^{-1}\left(\sup_{w \in A} S(w) - a\right),$$

$$\mathbb{E}|W - \mu|^3 \leq \lambda^{-1}(8q_{\max} + 1 + \sigma + \mathbb{E}\{|R|(W - \mu)^2\}).$$

PROOF.  The estimates for the variance are immediate from equality (3.23) and the bounds

$$\inf_{w \in A} S(w) \leq \mathbb{E}S(W) \leq \sup_{w \in A} S(w).$$

Note now that, from equation (3.22),

$$\mathbb{E}\{(W - \mu)g(W)\} = \lambda^{-1}\mathbb{E}\{S(W)\Delta g(W)\} + \lambda^{-1}\mathbb{E}\{Rg(W)\}$$

for all functions $g$, for which the expectations exist. With $K_\mu(w) = I[w > \mu] - I[w \leq \mu]$ and $g(w) = K_\mu(w)(w - \mu)^2$, we thus obtain

$$\mathbb{E}|W - \mu|^3 = \lambda^{-1}\mathbb{E}\{S(W)[(W - \mu)^2 + 2(W - \mu) + 1]\Delta K_\mu(W)\}$$
$$+ \lambda^{-1}\mathbb{E}\{S(W)(2(W - \mu) + 1)K_\mu(W)\}$$
$$+ \lambda^{-1}\mathbb{E}\{R(W - \mu)^2 K_\mu(W)\} =: B_1' + B_2' + B_3'.$$

Note now that $|K(w)| = 1$ and

$$\Delta K_\mu(w) = \begin{cases} 2, & \text{if } w = \lfloor\mu\rfloor, \\ 0, & \text{else,} \end{cases}$$

and thus, as $|\lfloor\mu\rfloor - \mu| \leq 1$ and $|S(w)| \leq 1$,

$$|B_1'| \leq 8\lambda^{-1}q_{\max},$$
$$|B_2'| \leq \lambda^{-1} + \lambda^{-1}\sigma.$$

The bound on $B_3'$ is immediate.  $\square$

**4. Applications.** In this section we illustrate our results using some examples in which $W = \sum_{i=1}^{n} J_i$ for a sequence $J = (J_1, J_2, \ldots, J_n)$ of random indicators. In [4] and [15], cases are considered where the $J_i$ have a local dependence structure; in contrast, the examples in this paper exhibit global dependence.

For latter use, we recall the following easy to prove fact.

LEMMA 4.1.    *Let $f : \mathbb{R} \to \mathbb{R}$ be a Lipschitz continuous function with Lipschitz constant $L_f$. Then, for any random variable $X$,*

$$\operatorname{Var} f(X) \le L_f^2 \operatorname{Var} X.$$

4.1. *Hypergeometric distribution.*    Assume that we have $N$ urns and $m$ balls, and that we distribute the balls uniformly into the $N$ urns, in such a way that there is at most one ball per urn. Clearly, the number of balls $W$ in the first $n$ urns has the hypergeometric distribution $\operatorname{Hyp}(m, n, N)$, for which

$$\sigma^2 = \operatorname{Var} W = \frac{nm(N-n)(N-m)}{(N-1)N^2}.$$

THEOREM 4.2.    *If $W$ has the hypergeometric distribution, then (3.1) and (3.2) hold with $R \equiv 0$ and $\lambda = \frac{N}{m(N-m+1)}$ and we have*

$$(4.1) \qquad \operatorname{Var} S \le \frac{nm(m+n)^2(N-n)(N-m)}{m^2(N-m+1)^2(N-1)N^2}.$$

*Thus, with $N = N(n) \asymp n$ and $m = m(n) \asymp n$,*

$$d_{\mathrm{TV}}(\mathscr{L}(W), \mathrm{TP}(\mu, \sigma^2)) = O(n^{-1/2}),$$

$$d_{\mathrm{loc}}(\mathscr{L}(W), \mathrm{TP}(\mu, \sigma^2)) = O(n^{-1}).$$

PROOF.    Consider the following construction. Pick uniformly an urn with a ball, and put this ball into any empty urn (including the urn from which the ball was picked). Denote now by $W'$ the number of balls in the first $n$ urns. Exchangeability of $(W, W')$ is easy to see and condition (2.2) is clearly satisfied. Now, $W' - W = 1$ is the event that a ball is picked from one of the urns $n+1, \ldots, N$ and put into one of the empty urns $1, \ldots, n$, thus,

$$\begin{aligned}
S(W) &= \mathbb{P}[W' = W + 1 | W] \\
(4.2) \qquad &= \frac{m - W}{m} \times \frac{n - W}{N - m + 1} \\
&= \frac{mn - (m+n)W + W^2}{m(N - m + 1)}
\end{aligned}$$

and, conversely,

$$\mathbb{P}[W' = W - 1 | W] = \frac{W}{m} \times \frac{N - n - m + W}{N - m + 1},$$

thus,

$$\mathbb{E}^W(W' - W) = \mathbb{E}^W I[W' - W = 1] - \mathbb{E}^W I[W' - W = -1]$$

$$= \frac{mn - NW}{m(N - m + 1)},$$

and (2.1) is satisfied with $R \equiv 0$ and $\lambda = \frac{N}{m(N - m + 1)}$.

Note now from (4.2) that $S$, as a function of $W$, is Lipschitz continuous with constant $L_S = \frac{m + n}{m(N - m + 1)}$; thus, applying Lemma 4.1, we have

$$\operatorname{Var} S \leq \frac{(m + n)\sigma^2}{m^2(N - m + 1)^2}. \qquad \square$$

This is enough to prove the $d_{\mathrm{TV}}$-order and, together with Corollary 3.9, the order $O(n^{-3/4})$ for the $d_{\mathrm{loc}}$-metric. Now, noting that Lemma 3.13 yields $\mathbb{E}|W - \mu|^3 = O(n^{3/2})$, we obtain from Theorem 3.11 the desired order $O(n^{-1})$ for the $d_{\mathrm{loc}}$-metric.

### 4.2. A parity problem.

Let $J_1, \ldots, J_n$ be a sequence of independent $\mathrm{Be}(1/2)$-distributed random indicators. Define

$$J_{n+1} := \begin{cases} 1, & \text{if } \sum_{i=1}^{n} J_i \text{ is odd,} \\ 0, & \text{else,} \end{cases}$$

and $V := \sum_{i=1}^{n+1} J_i$, so $V$ is simply obtained by "rounding" a $\mathrm{Bi}(n, 1/2)$-distributed random variable to the next even integer. An approximation of $V$ by a translated Poisson distribution will clearly not succeed; however, we may try with $W := \frac{1}{2} V$.

Regard now the following exchangeable pair coupling. Pick two random indices $K, L$ uniformly on $\{1, \ldots, n+1\}$ so that almost surely $K \neq L$, and define

$$(4.3) \qquad V' = V + 2 - 2J_K - 2J_L;$$

that is, take two summands of $V$ at random, and replace each of them by its complement.

LEMMA 4.3. *The pair $(V, V')$ defined as above is an exchangeable pair and $(W, W') := (\frac{1}{2} V, \frac{1}{2} V')$ satisfies (2.1) and (2.2) with $\lambda = 2/(n+1)$.*

Proof. It is enough to regard the situation on $M = \{0, 1\}^n$ because the values $J_1, \ldots, J_n$ uniquely determine the random variable $J_{n+1}$. Note first that construction (4.3) gives rise to a discrete time Markov chain on $M$, with jumps from $j \in M$ to $j' \in M$, if $j'$ differs from $j$ in exactly one or two coordinates ($j'$ differing in exactly one coordinate corresponds to $K$ or $L$ being equal to $n+1$). Now, as the jump from $j$ to $j'$ happens with the same probability as from $j'$ to $j$ and all the states are connected, it is easy to see that the such defined Markov chain is irreducible and reversible and that the equilibrium distribution assigns equal probability to any $j \in M$, which corresponds to $n$ independent $\text{Be}(1/2)$ random variables. Thus, exchangeability is proved.

Note now that

$$\mathbb{E}^J(V' - V) = 2 - \frac{2}{n(n+1)} \sum_{k=1}^{n+1} \sum_{\substack{l=1 \\ l \neq k}}^{n+1} (J_k + J_l)$$

$$= 2 - \frac{2}{n(n+1)} 2nV = 2 - \frac{4V}{n+1},$$

thus, we can take $\lambda = 2/(n+1)$. □

THEOREM 4.4. *For $W$ defined as above, (3.1) and (3.2) hold with $R \equiv 0$ and $\lambda = 2/(n+1)$ and if $n \geq 2$, we have $\sigma^2 = (n+1)/16$ and*

$$\text{Var } S \leq \frac{(4n-2)^2(n+1)}{16n^2(n+1)^2};$$

*thus, as $n \to \infty$,*

$$d_{\text{TV}}(\mathscr{L}(W), \text{TP}(\mu, \sigma^2)) = O(n^{-1/2}),$$
$$d_{\text{loc}}(\mathscr{L}(W), \text{TP}(\mu, \sigma^2)) = O(n^{-1}).$$

Proof. Note first that if $n \geq 2$, the $J_i$ are uncorrelated and, thus,

$$\sigma^2 = \text{Var}(V)/4 = (n+1)/16.$$

Now,

$$\mathbb{E}^J I[W' - W = 1] = \frac{1}{n(n+1)} \sum_{k=1}^{n+1} \sum_{\substack{l=1 \\ l \neq k}}^{n+1} (1 - J_k)(1 - J_l)$$

$$= \frac{n(n+1) - (4n-2)W + 4W^2}{n(n+1)} =: S(W).$$

Observe that $S$, as a function of $W$, is Lipschitz continuous with $L_S = \frac{4n-2}{n(n+1)}$; thus, applying Lemma 4.1,

$$\operatorname{Var} S(W) \leq \frac{(4n-2)^2 \sigma^2}{n^2(n+1)^2}.$$

This is enough to prove the $d_{\mathrm{TV}}$-order and, together with Corollary 3.9, the order $O(n^{-3/4})$ for $d_{\mathrm{loc}}$. Now, noting that Lemma 3.13 yields $\mathbb{E}|W - \mu|^3 = O(n^{3/2})$, we obtain from Theorem 3.11 the desired order $O(n^{-1})$ for the $d_{\mathrm{loc}}$-metric. $\quad\square$

4.3. *Anti-voter model on finite graphs.* We closely follow the setup of [14]; see also references therein and [13]. Let $G$ be an $n$-vertex $r$-regular graph, which is neither bipartite nor an $n$-cycle. At each vertex $i$ we assume that there is a "voter" attached, having an opinion $J_i^{(t)}$ which can take the values 0 or 1 in every time point $t \in \mathbb{N}$. Define a Markov chain by the following transition rule. Choose uniformly a random vertex, say, $i$; then, out of the neighborhood $\mathcal{N}_i$ of $i$, choose uniformly a random vertex, say, $j$, and let $J_i^{(t+1)}$ be the opposite of $J_j^{(t)}$ and leave the other voters untouched. Assume now that the Markov chain is in its equilibrium and put $W = \sum_{i=1}^n J_i := \sum_{i=1}^n J_i^{(0)}$.

THEOREM 4.5. *For the anti-voter model as described above, (3.1) and (3.2) hold with $R \equiv 0$ and $\lambda = 2/n$ and we have*

$$(4.4) \qquad \operatorname{Var} S \leq \frac{16r^2\sigma^2 + \operatorname{Var} Q}{16r^2n^2},$$

*where*

$$Q = \sum_{i=1}^n \sum_{j \in \mathcal{N}_i} (2J_i - 1)(2J_j - 1);$$

*hence, as $n \to \infty$,*

$$d_{\mathrm{TV}}(\mathscr{L}(W), \mathrm{TP}(\mu, \sigma^2)) = O\left(\frac{\sqrt{\operatorname{Var} Q}}{r\sigma^2} + \frac{1}{\sigma}\right),$$

$$d_{\mathrm{loc}}(\mathscr{L}(W), \mathrm{TP}(\mu, \sigma^2)) = O\left(\frac{(\operatorname{Var} Q)^{3/4}}{r^{3/2}\sigma^3} + \frac{1}{\sigma^{3/2}}\right).$$

REMARK 4.6. Note that the bound for $d_{\mathrm{TV}}$ in Theorem 4.5 is very similar to the bound for the weaker Kolmogorov metric $d_{\mathrm{K}}$ given in Theorem 1.3 of [14]; they obtain

$$(4.5) \qquad d_{\mathrm{K}}(\mathscr{L}(W_c), \mathcal{N}(0,1)) = O\left(\frac{\sqrt{\operatorname{Var} Q}}{r\sigma^2} + \frac{n}{\sigma^3}\right),$$

*where $W_c = (W - \mu)/\sigma$.*

EXAMPLE 4.7. Consider the sequence $K_n$ of complete graphs of size $n$. Rinott and Rotan [14] show that $\sigma^2 \asymp n$ and $\operatorname{Var} Q = O(n^3)$. Thus, from Theorem 3.1, the $d_{\mathrm{TV}}$-distance is of the order $O(n^{-1/2})$ and the $d_{\mathrm{loc}}$-distance of order $O(n^{-3/4})$ which proves the LLT. Now, from (4.8) below,

$$(4.6) \qquad S^*(J) = \frac{n(n-1) - (2n-1)W + W^2}{n(n-1)} = S(W),$$

and we can thus take $L_S = \frac{2}{n-1}$. From Lemma 3.13, we obtain $\mathbb{E}|W - \mu|^3 = O(n^{3/2})$ and, therefore, Theorem 3.11 yields the order $O(n^{-1})$ for $d_{\mathrm{loc}}$. Note that the estimates on $L_S$ are obtained only because of the explicit representation (4.6); they are difficult to obtain in general. For further examples of graphs, see [14].

PROOF OF THEOREM 4.5. Define $W' := \sum_{i=1}^n J_i^{(1)}$, and note that $(W, W')$ is an exchangeable pair, satisfying (2.1) and (2.2) with the choices $\lambda = 2/n$ and $R \equiv 0$ (for more details, see [14]). Now, let $K$ be the random index of the vertex that was resampled in the transition from $W$ to $W'$. As $W' = W - J_K + J_K^{(1)}$,

$$
\begin{aligned}
S^*(J) &= \mathbb{E}^J I[W' - W = 1] \\
&= \frac{1}{n} \sum_{i=1}^n \mathbb{E}^J \{ I[J_i = 0, J_i^{(1)} = 1] | K = i \} \\
&= \frac{1}{n} \sum_{i=1}^n (1 - J_i) \mathbb{E}^J \{ J_i^{(1)} | K = i \} \\
&= \frac{1}{n} \sum_{i=1}^n (1 - J_i) \left( 1 - \frac{1}{r} \sum_{j \in \mathcal{N}_i} J_j \right).
\end{aligned}
$$
(4.7)

With $X_i = 2J_i - 1$ and $\tilde{W} = \sum_{i=1}^n X_i$, (4.7) becomes

$$
\begin{aligned}
S^*(J) &= \frac{1}{4rn} \sum_{i=1}^n (1 - X_i) \left( r - \sum_{j \in \mathcal{N}_i} X_j \right) \\
&= \frac{1}{4rn} \left( rn - \sum_{i=1}^n \sum_{j \in \mathcal{N}_i} X_j - r \sum_{i=1}^n X_i + Q \right) \\
&= \frac{rn - 2r\tilde{W} + Q}{4rn}.
\end{aligned}
$$
(4.8)

The variance of $S^*$ is thus

$$\operatorname{Var} S^*(J) = \frac{\operatorname{Var}(2r\tilde{W}) + \operatorname{Var} Q - 4r \operatorname{Cov}(\tilde{W}, Q)}{16r^2 n^2}$$

(4.9)
$$= \frac{16r^2\sigma^2 + \operatorname{Var} Q}{16r^2 n^2},$$

because $\mathbb{E}\{X_i X_j X_k\} = 0$ for any choice of $i$, $j$ and $k$, due to the symmetry of the anti-voter model, and, hence, $\mathbb{E}\{\tilde{W}Q\} = 0$. $\quad\square$

**Acknowledgments.** I thank Andrew Barbour and the referee for many helpful suggestions and comments.

INSTITUT FÜR MATHEMATIK
UNIVERSITÄT ZÜRICH
WINTERTHURERSTRASSE 190
8057 ZÜRICH
SWITZERLAND
E-MAIL: adrian.roellin@math.unizh.ch